\documentclass[a4,leqno, 12pt]{article}
\usepackage{amsmath}
\newcommand{\N}{{\rm I\kern-.22em N}}
\newcommand{\Z}{{\sf Z\kern-.42em Z}}
\newcommand{\R}{{\rm I\kern-.22em R}}
\newcommand{\BbbC}
  {{\rm\kern.22em\rule[.1ex]{.06em}
       {1.4ex}\kern-.28em C}}
\newcommand{\BbbQ}
  {{\rm\kern.22em\rule[.1ex]{.06em}
       {1.4ex}\kern-.28em Q}}

\newcommand\qed{$\Box$}

\newtheorem{thm}{Theorem}[section]
\newtheorem{lemma}[thm]{Lemma}
\newtheorem{prop}[thm]{Proposition}

\makeatletter

\@addtoreset{equation}{section}
\makeatother

\date{}
\title
{
Finite speed of propagation 
     in 1-D degenerate \\ Keller-Segel system
}

\author{\vspace*{5mm} 
\\
{\small 
Department of Mathematics,
       Tsuda University,} \\
{\small       2-1-1, Tsuda-chou, Kodaira-shi, Tokyo, 187-8577, {\sc Japan}},
\\
{\small sugiyama@tsuda.ac.jp}
       }
\begin{document}
\maketitle
 \vskip5mm
\begin{abstract}
We consider the following Keller-Segel system of degenerate type:
\begin{equation}\label{eqn;KS} 
\tag{KS}  \left\{
    \begin{array}{llll}
         & \displaystyle \frac{\partial u}{\partial t} = 
             \frac{\partial}{\partial x} 
             \Big( \frac{\partial u^m}{\partial x} 
                   - u^{q-1} \cdot \frac{\partial v}{\partial x} \Big),  
                    & x \in \R, \ t>0, \nonumber \\
         & \displaystyle 0 = 
	  \frac{\partial^2 v}{\partial x^2} - \gamma v + u, 
                    & x \in \R, \ t>0, \nonumber \\
         & u(x,0)   = u_0(x),   \quad 
                    & x \in \R, 
     \end{array}
  \right.
\end{equation}
where $m>1, \ \gamma > 0, \ q \ge 2m$. 
We shall first construct a weak solution $u(x,t)$ of (KS) such that 
$u^{m-1}$ is Lipschitz continuous and such that
$\displaystyle u^{m-1+\delta}$
for $\delta>0$ is of class $C^1$ with respect to the space variable $x$.
As a by-product, we prove
the property of finite speed of propagation 
   of a weak solution $u(x,t)$ of (KS), {\it i.e.,}  
that 
a weak solution $u(x,t)$ of (KS) has a compact support in 
$x$ for all $t>0$ 
if the initial data $u_0(x)$ has a compact support in $\R$.   
We also give 
both upper and lower bounds of the interface of the weak solution $u$ of (KS).
\end{abstract}
\vspace{2mm}

%
%   section  1
%

\section{Introduction}
\qquad
We consider the following Keller-Segel system of degenerate type:
\begin{equation}\label{eqn;KS} 
\tag{KS}  \left\{
    \begin{array}{llll} 
         & \displaystyle \frac{\partial u}{\partial t} = 
             \frac{\partial}{\partial x} 
             \Big( \frac{\partial u^m}{\partial x} 
                   - u^{q-1} \cdot \frac{\partial v}{\partial x} \Big),  
                    & x \in \R, \ t>0, \nonumber \\
         & \displaystyle 
0 = 
	  \frac{\partial^2 v}{\partial x^2} - \gamma v + u, 
                    & x \in \R, \ t>0, \nonumber \\
         & u(x,0)   = u_0(x),   \quad 
                    & x \in \R, 
     \end{array}
  \right.
\end{equation}
where $m>1, \ \gamma > 0, \ q \ge 2m$.
The initial data $u_0$ is a non-negative function and 
in $L^{1} \cap L^{\infty}(\R)$ with 
              $u_0^m \in H^1(\R). $
This equation is often called as the Keller-Segel model describing 
the motion of the chemotaxis molds. (see {\it e.g.,} \cite{ChPe}.)

The aim of this paper is to construct a weak solution $u(x,t)$ of (KS) such that 
$u^{m-1}$ is Lipschitz continuous 
and such that $\displaystyle u^{m-1+\delta}$
for $\delta>0$
is of class $C^1$ with respect to the space variable $x$.
The regularity property 
whether $u^{m-1}$ is Lipschitz continuous or of class $C^1$ 
plays an important role for the investigation of the behaviour of the interface
to the solution $u$ of (KS).
Our result shows that the power $m-1$ to $u$ exhibits the borderline behaviour
between Lipschitz continuity and $C^1$-regularity. 
Indeed, as a by-product of Lipschitz continuity for $u^{m-1}$, 
we prove that
a weak solution $u(x,t)$ of (KS) possesses 
the property of finite speed of propagation {\it i.e.,}  that 
a weak solution $u(x,t)$ of (KS) has a compact support in 
$x$ for all $t>0$ 
if the initial data $u_0(x)$ has a compact support in $\R$. 

Similar results have been obtained for the porous medium equation:
\begin{equation*}
\tag{PME} 
  \left\{
    \begin{array}{llll} 
\displaystyle 
\frac{\partial U}{\partial t} 
       & = &
\displaystyle
             \frac{\partial^2 U^m}{\partial x^2},
                 \hspace{2cm}  &  x \in \R, \ t>0,
\nonumber \\
U(x,0) & = & U_0(x), 
                 \hspace{2cm}  &  x \in \R.
     \end{array}
  \right.
\end{equation*}
It is known that the comparison principle gives both upper and lower
bounds of all solutions $U$ to (PME) by means of the {\it Barenblatt solution}
$V_B$ which is an exact solution of (PME).
Hence the property of finite speed of propagation of $U$ is a direct consequence of
the explicit form of $V_B$ since supp $V_B(\cdot,t)$ is compact in $\R$ for all
time $t$.

Our purpose is to prove the property of finite speed of propagation
for (KS) to which the comparison principle is not available.
To this end, one makes use of the notion of the domain of dependence  
which is useful for the proof of uniqueness 
                       of solutions to the linear wave equations. 
For instance, 
the half-cone like region $D_T$ defined by
\begin{eqnarray*}
D_T 
& := & 
\Big\{ (x,t) 
; \ 
        -c t +a \le x \le c t +b, \quad 0 \le t <T
          \Big\},
\quad a<b, \ c>0
\end{eqnarray*}
makes it possible to prove that the solution of the linear wave equation
with the propagation speed $c$
vanishes on $D_T$ for the initial data $u_0$ such that $u_0(x) \equiv 0$
on $I \equiv [a,b]$.

To deal with (KS), we generalize such an idea, and consider
the {\it curved} half-cone like region. 
Indeed, suppose that $u_0(x)=0$ on $I$.
Then our curved half-cone like region $D_T$ with respect to $I$ can be
expressed by 
\begin{eqnarray}\label{region-DT}
D_T 
& := & 
\Big\{ (x,t) 
; \ 
        \xi(t) \le x \le \Xi(t), \quad 0 \le t <T
          \Big\},
\end{eqnarray}
where $\xi(t)$ and $\Xi(t)$ are the solutions of 
the following initial value problems:
\begin{equation} 
\tag{IE} \left\{
\begin{array}{lll}
\xi^{'}(t) 
& = & - \frac{\partial}{\partial x}
             \Big(\frac{m}{m-1} u^{m-1} \Big)(\xi(t),t)
      + u^{q-2}
              \cdot \frac{\partial v}{\partial x} (\xi(t),t), \quad \xi(0)=a,
\label{int;1} \\
\Xi^{'}(t) 
& = & - \frac{\partial}{\partial x}
             \Big(\frac{m}{m-1} u^{m-1} \Big)(\Xi(t),t)
      + u^{q-2}
              \cdot \frac{\partial v}{\partial x} (\Xi(t),t), \quad \Xi(0)=b.
\end{array}
\right.
\end{equation}
Unfortunately, Lipschitz continuity of $u^{m-1}$ is too weak to ensure
the existence of solutions $\{\xi(t), \Xi(t)\}$ to (IE).
Hence we need to regularize $u$ by $u_{\varepsilon}$ with small
parameter $\varepsilon>0$, 
and deal with the approximating solutions 
       $\{\xi_{\varepsilon}(t), \Xi_{\varepsilon}(t)\}$ which correspond
to (IE) with $u$ replaced by $u_{\varepsilon}$.
It is shown that Lipschitz continuity of $u^{m-1}$ guarantees the
existence of uniform limit $\{\xi(t), \Xi(t)\}$ on $0 \le t \le T$ of  
    $\{\xi_{\varepsilon}(t), \Xi_{\varepsilon}(t)\}$ as $\varepsilon \to 0$.
Then we see that $u(x,t)=0$ on $D_T$.
\vspace{6mm} \\
\qquad
Our definition of a weak solution to (KS) now reads:
\vspace{2mm} \\
{\bf Definition 1} \ 
{\it  Let $m,\gamma$ and $q$ be constants as
$m>1, \ \gamma> 0, \ q \ge 2. 
$ 
Let $u_0$ be a non-negative function in $\R$ with
$u_0 \in L^{1} \cap L^{\infty}(\R)$ 
and $u_0^m \in H^1(\R)$.
 A pair of non-negative functions  
$(u,v)$ defined in $\R \times [0,T)$ 
is said to be a weak solution of {\rm (KS)} on $[0,T)$
if
\vspace{2mm} \\
{\rm i)} \ \ \ 
$u \in  L^{\infty} (0,T; L^2(\R)),$ 
$u^m \in L^2(0,T;H^1(\R)) 
,$  
\vspace{3mm} \\
{\rm ii)} \ \ 
$
v \in  L^{\infty}(0,T; H^2(\R)),
$
\vspace{3mm} \\
{\rm iii)} \  $(u,v)$ satisfies (KS)
         in the sense of distributions: 
                     {\it i.e.,}
\begin{eqnarray*}
&& \int_0^T \int_{\R} 
  \left(\partial_x  u^m \cdot \partial_x  \varphi 
        - u^{q-1} \partial_x  v \cdot \partial_x  \varphi
        - u \cdot \partial_t \varphi \right) \ dxdt 
  \ = \  \int_{\R} u_0(x) \cdot \varphi (x,0) \ dx, \\
&& 
\hspace{6.5cm}
{\it for \ all \ functions}  \ \varphi \in C_0^{\infty}(\R \times [0,T)),
\\
&& \\
&&  
      -\partial_x^2  v  + \gamma v  - u =0 
\qquad {\it for \ a.a.}  \ (x,t) \ {\it in } \ \R \times (0,T).
\end{eqnarray*}
}
\vspace{2mm} \\
\qquad
Concerning the local-in-time existence of weak solutions to (KS),
the following result can be shown by a slight modification of argument 
    developed
           by the author \cite[Theorem 1.1]{Sblow-uq-2005}.
\begin{prop} 
[{\small local existence of weak solution and 
     its $L^{\infty}$ uniform bound}
]
\label{pr;local existence} 
\ \\
Let $m>1,\ \gamma > 0, \ q \ge 2$.
Suppose that the initial data $u_0$ is non-negative everywhere.
Then, 
 ${\rm (KS)}$ 
       has a non-negative weak solution 
         {\rm ($u, v$)}
on $[0,T_0)$ with
$
T_0 = \Big( \| u_{0} \|_{L^{\infty}(\R)} + 2 \Big)^{-q}. 
$  
Moreover, $u(t)$ satisfies the following {\it a priori} estimate
\begin{eqnarray} \label{by}
\| u(t) \|_{L^{\infty}(\R)} 
   & \le & \| u_0 \|_{L^{\infty}(\R)} + 2  
\ \   \quad {\rm for \ all} \ t \in [0,T_0).  
\end{eqnarray}
\end{prop}
\ \vspace{2mm} \\
{\bf Remark 1.} \
Concerning the global-in-time existence of weak solutions to (KS),
the author and Kunii \cite{KS-SuKu} obtained the following result:
Let $m,\gamma,q$ and the initial data $u_0$ be as in Definition 1.
In the case $q<m+2$,
there exists a weak solution $u$ of (KS) on $[0,\infty).$
On the other hand, in the case $q \ge m+2$,
the weak solution $u$ of (KS) on $[0,\infty)$ can be constructed provided
$\|u_0\|_{L^{\frac{q-m}{2}}(\R)}$ is sufficiently small. 
\vspace{4mm} \\
\qquad
Now, we construct a weak solution $u$ of (KS) 
   with some additional regularity for the velocity potential $u^{m-1}$.
\begin{thm} \label{thm:lip}
Let $m>1, \gamma > 0$ and $q \ge 2m$.
Let the initial data $u_0$ be as in Definition 1. 
In addition, we assume that $u_0^{m-1}$ is Lipschitz continuous in $\R$.
Then, the weak solution $u$ of (KS) on $[0,T_0)$ 
given by Proposition \ref{pr;local existence}
has the following additional properties (i) and (ii): 
\vspace{2mm} \\
(i) \ $u^{m-1}(x,t)$ is Lipschitz continuous with respect to $x$ for all $0 \le t < T_0$ 
with the estimate
\begin{eqnarray}\label{pr-adj}
\sup_{0<t<T_0}\| \partial_x u^{m-1}(t) \|_{L^{\infty}(\R)}
& \le & C,
\end{eqnarray}
where $C=C(m,\gamma,q,u_0)$.
\vspace{2mm} \\
(ii) \ For every $\delta>0$, 
$u^{m-1+\delta}(\cdot,t)$ is a $C^1$-function with respect to $x$ 
for all $0<t<T_0$ with the property that
$\partial_x u^{m-1+\delta}(x,t)=0$ at the points $(x,t) \in \R \times
 (0,T_0)$ such that $u(x,t)=0$. 
Furthermore, in the case of $1<m<2$, we have
$\partial_x u(x,t)=0$ at the same points $(x,t)$ as above. 
\end{thm} 
{\bf Remark 2.}
\ 
(i) \ By the fundamental inequality
\begin{eqnarray} \label{holder-est}
&& |u(x,t)-u(y,t)| \ \le \ 
\left\{
\begin{array}{lcl}
\frac{2^{1/m-1}}{m-1} 
\|u\|_{L^{\infty}(Q_{T_0})}^{2-m}|u^{m-1}(x,t)-u^{m-1}(y,t)|, & \ 1<m<2,& \\
|u^{m-1}(x,t)-u^{m-1}(y,t)|^{\frac{1}{m-1}}, & \ m \ge 2 &
\end{array}
\right.
\end{eqnarray} 
for all $x,y \in \R, \ 0<t<T_0,$ we have by Theorem \ref{thm:lip} that
for every $0<t<T_0$, 
$u(\cdot,t)$ is a H\"{o}lder continuous function in $\R$ with the exponent 
$\mu=\min\{1,\frac{1}{m-1}\}$.
\vspace{2mm} \\
(ii) \ For (PME), it is well-known that
$\partial_x U^{m-1}(\cdot,t)$ becomes a discontinuous
function in $\R$ after some definite time $t$.
Our result in Theorem \ref{thm:lip} makes it clear that continuity in
$x$ of $\partial_x u^p(x,t)$ is guaranteed for all $p>m-1$ and all $0<t<T_0$.
It seems to be an interesting question whether $\partial_x
u^{m-1}(\cdot,t)$ is really discontinuous in $\R$ or not.
\vspace{2mm} \\
(iii) \ 
The hypothesis $q \ge 2m$ seems to be redundant.
Indeed, such restriction on $q$ stems from 
choice of the transformation $\psi$ in (\ref{PSIAT})
which may have a certain freedom to apply the Bernstein method 
to the uniform estimate of $\partial_x u^{m-1}$ for $\varepsilon$.
It should be noted that in the case $\gamma=0$, 
we can relax this restriction to $q \ge m+1$. See (\ref{JJ}) below.
\vspace{2mm} \\
\vspace{2mm} \\
\qquad
By (\ref{pr-adj}) in Theorem \ref{thm:lip}, 
we can construct a pair of continuous functions $\xi(t)$ and $\Xi(t)$
on $[0,T_0)$ such that the region  $D_{T_0}$ defined by (\ref{region-DT})
belongs to the interior of the domain surrounded by the interface of $u$, 
which leads us to the property of finite speed of propagation to (KS).
\begin{thm}[property of finite propagation speed]
\label{thm;finite speed}
Let $m>1, \ \gamma > 0$ 
and $q \ge 2m$.
Let the initial data $u_0$ be as in Definition 1. 
In addition, we assume that 
$u_0(x)=0$ on some interval $I \equiv [a,b]$
and that $u_0^{m-1}$ is Lipschitz continuous in $\R$.
Suppose that $u$ is the weak solution of (KS) on $[0,T_0)$ given by
 Theorem \ref{thm:lip}.
Then, there exists a pair $\{\xi(t),\Xi(t)\}$ of continuous functions 
on $[0,T_0)$ with the following properties (i) and (ii):
\vspace{2mm} \\
(i) \ $\xi, \Xi \in W^{1,\infty}(0,T_0)$ \ with $\xi(0)=a, \ \Xi(0)=b$;
\vspace{2mm} \\
(ii) \ $u(x,t)=0$ \ \ for \ $\xi(t) \le x \le \Xi(t), \ 0 \le t<T_0$.
\end{thm}
{\bf Remark 3.}
\
(i) \
Concerning (PME), 
the interface of $U$ 
    can be explicitly determined by the solutions 
             $\hat{\xi}(t)$ and $\hat{\Xi}(t)$ of
                             the following initial value problems:
\begin{equation} 
\left\{
\begin{array}{lll}
\hat{\xi}^{'}(t) 
& = & - \frac{\partial}{\partial x}
             \Big(\frac{m}{m-1} U^{m-1} \Big)(\hat{\xi}(t),t),
\quad
\hat{\xi}(0)=a,
\nonumber \\
\hat{\Xi}^{'}(t) 
& = & - \frac{\partial}{\partial x}
             \Big(\frac{m}{m-1} U^{m-1} \Big)(\hat{\Xi}(t),t),
\quad \hat{\Xi}(0)=b.
\nonumber 
\end{array}
\right.
\end{equation}
Indeed, by the comparison principle Knerr \cite{Kn} showed that 
if $U_0(x) = 0$ on some interval $I=[a,b]$ and $U_0(x)>0$ 
   on $I^c=\R \backslash I$, 
then it holds that
$U(x,t) = 0$ for $\hat{\xi}(t) \le x \le \hat{\Xi}(t)$ 
and $U(x,t)>0$ for $x < \hat{\xi}(t)$ and $x > \hat{\Xi}(t)$
 for all $0 \le t <\infty$.
We call such $\hat{\xi}(t)$ and $\hat{\Xi}(t)$ 
the {\it interface} of (PME).
\vspace{2mm} \\
(ii) \ 
Compared with (PME), it is not clear 
whether (IE) determines the exact interface of (KS) 
to which the comparison principle is not available. 
However, if 
$\xi_1(t)$ and $\Xi_1(t)$ are the interface of (KS), {\it i.e.,} that 
$\xi_1(t)$ and $\Xi_1(t)$ have the property that 
$$
u(x,t) = 0 {\rm \ \ in} \  I_{t} := [\xi_1(t),\Xi_1(t)] 
\quad {\rm and} \quad u(x,t)>0 \ {\rm in  \ some \ neighbourhood \ outside \
of } \ I_{t} 
$$
for all $0 \le t <T_0$,
then we can see that $\xi(t)$ and $\Xi(t)$ given by Theorem \ref{thm;finite speed}
satisfy the estimates
$$
\xi_1(t) \le \xi(t), \qquad \Xi(t) \le \Xi_1(t) \quad {\rm for \ all} \
0 \le t < T_0.
$$
Hence our result may be regarded as an estimate of the maximum and the
minimum of the interface of (KS).
Other observations were done 
by Mimura-Nagai \cite{MN} and Bonami-Hilhorst-Logak-Mimura \cite{B-H-L-M-2}.
\vspace{4mm} \\
\qquad
This paper is organized as follows.
In Section 2, 
we shall first recall the approximating problem (KS)$_{\varepsilon}$ of (KS) 
                introduced by \cite{KS-SuKu}.
Our main purpose is devoted to the derivation of 
               uniform gradient bound with respect to $\varepsilon>0$ of 
   the approximating velocity potential $w_{\varepsilon}=\frac{m}{m-1} u_{\varepsilon}^{m-1}$, 
    where $u_{\varepsilon}$ is the smooth solution
                        of (KS)$_{\varepsilon}$.
Bernstein's method plays an important role to obtain our uniform estimate.
(see {\it e.g.,} \cite{Li}.)
Then in Section 3, 
   by the standard compactness argument, we shall prove the Lipschitz
continuity of the velocity potential $w=\frac{m}{m-1}u^{m-1}$ for the
weak solution $u$ of (KS).
It is expected that $\partial_x w(x,t)$ becomes a discontinuous
function in $x$ after some finite time $t$.
However, we shall show that 
for $p>m-1$, $\partial_x u^p(x,t)$ is, in fact,
a continuous function in $\R$ for all $t \in [0,T_0)$.
Section 4 is devoted to the construction of continuous curves 
$\xi(t)$ and  $\Xi(t)$ such that $u(x,t)=0$ on $D_{T_0}$ defined by (\ref{region-DT}),
which implies the property of the finite speed of propagation for (KS).
\vspace{4mm} \\
\qquad
We will use the simplified notations: 
\vspace{2mm} \\
1) \ $\partial_t = \frac{\partial}{\partial t}, \ $
   \ $\partial_x = \frac{\partial}{\partial x}, \ $
   \ $\partial^2_{x} = \frac{\partial^2}{\partial x^2}, \ $
   \ $\partial^3_{x} = \frac{\partial^3}{\partial x^3}, \ $
\\
2) \
$\| \cdot \|_{L^r} = \| \cdot \|_{L^r(\R)},  
             (1 \le r \le \infty), \
    \int \cdot \ dx:= \int_{\R} \cdot \ dx,
$
\vspace{1mm} \\
3) \ $Q_T:= \R \times (0,T)$, \ 
\vspace{2mm} \\
4) \ 
When the weak derivatives $\partial_x u, \partial_x^2 u$ 
                         and $\partial_t u$ 
          are in $L^p(Q_T)$ for some $p \ge 1$,
we say that $u \in W^{2,1}_p(Q_T)$, {\it i.e.}, 
\begin{eqnarray*}
W^{2,1}_p(Q_T) 
& := & \Big\{u \in L^p(0,T;W^{2,p}(\R)) \cap
                  W^{1,p}(0,T;L^p(\R)) ; \\
&&   \quad 
       \| u \|_{W^{2,1}_p(Q_T)}
            := \|u\|_{L^p(Q_T)}+\|\partial_x u\|_{L^p(Q_T)}
                       + \|\partial_x^2u\|_{L^p(Q_T)}
                        + \|\partial_t u\|_{L^p(Q_T)} < \infty
         \Big\} .
\end{eqnarray*}
\section{{\bf Approximating Problem}} \label{3}
\hspace{5mm}
\quad 
In order to justify the formal arguments, 
we introduce the following approximating equations of
(KS):
\label{sect; time local}
$$
{\rm (KS)_{\varepsilon}} 
~\left\{\begin{array}{rclll}
\begin{displaystyle} 
\partial_t u_{\varepsilon} \end{displaystyle}(x,t)
  & = & 
    \partial_x 
       \Big(\partial_x  
              (u_{\varepsilon} + \varepsilon)^m                    
           -  (u_{\varepsilon}+\varepsilon)^{q-2} u_{\varepsilon} \cdot 
                \partial_x  v_{\varepsilon} \Big) , 
  & (x,t) \in \R \times (0,T),  
                       \vspace{2mm} \\
0
  & = & 
    \partial_x^2  v_{\varepsilon} - \gamma v_{\varepsilon} 
           +  u_{\varepsilon}, 
  &  (x,t) \in \R \times (0,T),    
 \vspace{2mm} \\
u_{\varepsilon}(x,0)  
      & = & u_{0\varepsilon}(x), \ 
                  &  x \in \R, 
\end{array}
\right.
\hspace{3cm}
$$
where $\varepsilon>0$ is a positive parameter. 
\vspace{2mm} \\
\qquad
Let us introduce the following assumption on the initial data
                  $u_{0\varepsilon}$ with $\varepsilon>0$.
\noindent
\begin{description}
\item[(A.1)] 
$u_{0\varepsilon} \ge 0$ for all $x \in \R$ 
              and $u_{0\varepsilon} \in W^{2,p}(\R)$ with
\begin{eqnarray*}
\sup_{0<\varepsilon<1}
   \| u_{0\varepsilon}\|_{L^p(\R)} 
            \le \|u_0\|_{L^p(\R)}
&&
       {\rm for \ all} \ p \in [1,\infty],  
\nonumber \\
\|u_{0\varepsilon} - u_0 \|_{L^p(\R)} \to 0 \quad {\rm as}
\ \varepsilon \to 0
&&
       {\rm for \ all} \ p \in [1,\infty).  
\end{eqnarray*} 
\item[(A.2)]  
$u_{0\varepsilon} \in W^{1,2}(\R)$ 
\quad
with
\
$
\sup_{0<\varepsilon<1}
   \| \partial_x u_{0\varepsilon} \|_{L^2(\R)} 
            \le \|\partial_x u_0 \|_{L^2(\R)}.
$
\end{description}
\vspace{4mm}
{\bf Definition 2} \quad
We call ($u_{\varepsilon}, v_{\varepsilon}$) 
a {\it strong solution} of (KS)$_{\varepsilon}$ 
      if it belongs to $W^{2,1}_p \times W^{2,1}_p(Q_T)$ 
                       for some $p \ge 1$
       and (KS)$_{\varepsilon}$ 
             is satisfied almost everywhere.
\vspace{4mm} \\
\qquad
For the strong solution, 
we consider the case $p=3$ and introduce the space ${\bf W}(Q_T)$ defined by
\begin{eqnarray}\label{def-wqt}
{\bf W}(Q_T)
& := & 
                 W^{2,1}_{3}
                       \times 
                 W^{2,1}_{3}(Q_T).
\end{eqnarray}
\qquad 
In \cite{Sblow-uq-2005}--\cite{KS-SuKu}, 
the following proposition concerning the
existence of the strong solution was proved :
\begin{prop} \label{pr;apto}
({\small local existence of approximating solution})
Let $m \ge 1, \ \gamma > 0, \ q \ge 2$.
We take $T_0:=(\|u_0\|_{L^{\infty}(\R)}+2)^{-q}$.
Then, for every $\varepsilon>0$ and every initial data
 $u_{0\varepsilon}$ satisfying the hypothesis {\rm (A.1)},
${\rm (KS)_{\varepsilon}}$ 
       has a unique non-negative strong solution 
         {\rm ($u_{\varepsilon}, v_{\varepsilon}$)} in ${\bf W}(Q_{T_0})$.
Moreover, $u_{\varepsilon}(t)$ satisfies 
the following {\it a priori} estimate
\begin{eqnarray} \label{estK0}
\| u_{\varepsilon}(t) \|_{L^{\infty}(\R)} 
   & \le & \| u_0 \|_{L^{\infty}(\R)} + 2 
 \   \quad {\rm for \ all} \ t \in [0,T_0) \ 
              {\rm and \ all} \ \varepsilon \in (0,1].
\end{eqnarray}
\end{prop}
{\bf Remark 4.} \
(i) \ It should be noted that the time interval $[0,T_0)$ of the
existence of the strong solution $(u_{\varepsilon},v_{\varepsilon})$ can
be taken uniformly with respect to $\varepsilon>0$.
\vspace{2mm} \\
(ii) \ The weak solution $(u,v)$ of (KS) on $[0,T_0)$ 
       given by Proposition \ref{pr;local existence} can be constructed
as the weak limit of $(u_{\varepsilon},v_{\varepsilon})$ 
as $\varepsilon \to 0$, where $(u_{\varepsilon},v_{\varepsilon})$ is the
strong solution in Proposition \ref{pr;apto}.
More precisely, 
by choosing a subsequence of $(u_{\varepsilon},v_{\varepsilon})$
which we denote by $(u_{\varepsilon},v_{\varepsilon})$ itself for
simplicity, we have
\begin{eqnarray*}
u_{\varepsilon} 
  \ \  \rightharpoonup \ \ u \qquad && {\rm weakly-star} 
    \ {\rm in} \ L^{\infty}(0,T_0; L^2(\R)),
\nonumber \\
u_{\varepsilon}^m 
  \ \ \rightarrow \ \ u^m \qquad && {\rm weakly } \ {\rm in} \ 
                L^2(0,T_0; H^1(\R)) \ {\rm and \ strongly \ in} \ 
                               C([0,T_0); L_{loc}^2(\R)),
\nonumber \\
v_{\varepsilon} 
   \ \ \rightharpoonup \ \ v \qquad && 
          {\rm weakly-star} \ {\rm in} \ 
                     L^{\infty}(0,T_0; H^2(\R))
\end{eqnarray*}
as $\varepsilon \to 0.$
In what follows, we assume that the sequence of approximating solutions
($u_{\varepsilon},v_{\varepsilon}$) satisfies the above convergence.
\vspace{2mm} \\
(iii) \ 
The strong solution $(u_{\varepsilon},v_{\varepsilon}) \in {\bf
W}(Q_{T_0})$ is more regular. Indeed, for every $\varepsilon>0$,
it can be shown that
$
u_{\varepsilon}, v_{\varepsilon} \in C^{\infty}(\R \times (0,T_0)).
$
\vspace{6mm} \\
\qquad
The following lemma gives the gradient estimate 
for the velocity potential $u^{m-1}$. 
\begin{lemma} \label{lip}
Let $m>1, \ \gamma > 0$ and $q \ge 2m$. 
Let the initial data $u_0$ be as in Definition 1. 
For every $\varepsilon>0$, we take 
$u_{0\varepsilon}$ so that the hypothesis {\rm (A.1)--(A.2)} are satisfied.
In addition, we assume that $u_{0\varepsilon}^{m-1}$ is Lipschitz continuous in $\R$.
Then the strong solution $u_{\varepsilon}$ of (KS)$_{\varepsilon}$ on $[0,T_0)$
given by Proposition \ref{pr;apto} has the following property:
\begin{eqnarray}\label{adj}
\sup_{0<\varepsilon<1}
\Big(
\sup_{0<t<T_0}\| \partial_x (u_{\varepsilon}+\varepsilon)^{m-1} \|_{L^{\infty}(\R)}
\Big)
& \le & C,
\end{eqnarray}
where $C=C(m,\gamma,q,u_0)$.
\end{lemma} 
{\it Proof of Lemma \ref{lip}.} \
For the sake of simplicity, we denote
$(u_{\varepsilon}, v_{\varepsilon})$ by $(u,v)$.
To treat the velocity potential, 
let us define $w := \frac{m}{m-1}(u+\varepsilon)^{m-1}$. 
Multiplying the first equation of (KS)$_{\varepsilon}$
                    by $m(u+\varepsilon)^{m-2}$
and then rewriting the resultant identity in terms of $w$, we have
\begin{eqnarray} \label{w}
\partial_t  w
& = &
        (m-1)w \cdot \partial_{x}^2 w
          + |\partial_x  w|^2
           - \Big( (q-2) (u+\varepsilon)^{q-3} u
               + (u + \varepsilon)^{q-2} 
             \Big) \cdot \partial_x  v 
                   \cdot \partial_x  w
\nonumber \\
&& -(m-1)(u+\varepsilon)^{q-3}u
                         \cdot \partial_x^2  v
                         \cdot w. 
\end{eqnarray}
Now we apply Bernstein's method. 
Introducing the convex transformation $\psi: \bar{w} \to w$, 
determined below (\ref{PSIAT}), 
we rewrite the identity (\ref{w}) by means of $\bar{w}=\psi^{-1}(w)$ in the
following form:
\begin{eqnarray} \label{apri nabla; bar u..}
\partial_t \bar{w}  
& = & (m-1) \psi \cdot 
       \Big(\frac{\psi^{''}}{\psi^{'}} |\partial_x \bar{w}|^2
                       + \partial_x^2 \bar{w} \Big)
            +
      \psi^{'} |\partial_x \bar{w}|^2
\nonumber \\
&& \quad       -
    \Big((q-2)(u+\varepsilon)^{q-3} u+ (u+\varepsilon)^{q-2} \Big)
                \cdot \partial_x v \cdot \partial_x \bar{w}
  - (m-1) \frac{\psi}{\psi^{'}} \cdot 
      (u+\varepsilon)^{q-3} u \cdot \partial_x^2 v.  
\end{eqnarray}
We note that 
\begin{eqnarray}
&& \psi(\bar{w}) 
\ = \ w \ = \ \frac{m}{m-1} (u+\varepsilon)^{m-1},
\quad
\psi^{'}(\bar{w}) \cdot \partial_x \bar{w}
\ = \ m(u+\varepsilon)^{m-2} \partial_x u,
\label{psi;w-1} \\
&& (u+\varepsilon)^{q-i} \cdot \partial_x u
\ = \ \frac{1}{m} \cdot \psi^{'} \cdot (u+\varepsilon)^{q-m-i+2} 
                            \cdot \partial_x \bar{w}
\quad {\rm for} \ i=3,4. 
\label{psi;w-2} 
\end{eqnarray}
Differentiating both sides of (\ref{apri nabla; bar u..})
   with respect to $x$,
we obtain from (\ref{psi;w-1}) and (\ref{psi;w-2}) that
\begin{eqnarray}\label{bar u 1..}
\partial_t \partial_{x} \bar{w}
& = & (m-1)
      \psi'(\bar{w}) \cdot 
           \Big(\frac{\psi^{''}(\bar{w})}{\psi^{'}(\bar{w})} 
                                |\partial_x \bar{w}|^2
                       + \partial_x^2 \bar{w} 
           \Big) \cdot \partial_x \bar{w}
                                    \nonumber \\
&& + \Big(
         (m-1) \psi(\bar{w}) \cdot   
    \Big(  \frac{ \psi''(\bar{w})}{\psi'(\bar{w})} \Big)^{'}
    + \psi^{''}(\bar{w})            
     \Big)
          \cdot (\partial_x  \bar{w})^3 
\nonumber \\
&&  + \ 2 \Big(
              (m-1) \psi(\bar{w}) \cdot \frac{\psi''(\bar{w})}{\psi'(\bar{w})}
                  + \psi^{'}(\bar{w})
          \Big)
                      \cdot \partial_x \bar{w} 
                          \cdot \partial_x^2  \bar{w} 
\  + \ (m-1) \psi(\bar{w}) 
                          \cdot \partial_x^3 \bar{w} 
                                     \nonumber \\
&& - \ (q-2)(q-3) \cdot \frac{1}{m} \cdot \psi^{'}(\bar{w})
         \cdot (u+\varepsilon)^{q-m-2} u 
         \cdot \partial_x v \cdot (\partial_x \bar{w})^2
\nonumber \\
&& - \ 2(q-2) \cdot \frac{1}{m} \cdot \psi^{'}(\bar{w})
         \cdot (u+\varepsilon)^{q-m-1} 
         \cdot \partial_x v \cdot (\partial_x \bar{w})^2
\nonumber \\
&& - \ (q-2) \cdot \frac{m-1}{m} \cdot \psi(\bar{w})
         \cdot (u+\varepsilon)^{q-m-2} u 
         \cdot \partial_x^2 v 
         \cdot \partial_x \bar{w} 
\nonumber \\
&& - \ (q-2) \cdot \frac{m-1}{m} \cdot \psi(\bar{w})
         \cdot (u+\varepsilon)^{q-m-2}  u
         \cdot \partial_x v \cdot 
                               \partial_x^2 \bar{w} 
\nonumber \\
&& - \   \frac{m-1}{m} \cdot \psi(\bar{w})
         \cdot (u+\varepsilon)^{q-m-1}  
         \cdot \partial_x^2 v \cdot \partial_x \bar{w} 
\nonumber \\
&& - \   \frac{m-1}{m} \cdot \psi(\bar{w})
         \cdot (u+\varepsilon)^{q-m-1}  
         \cdot \partial_x v  
         \cdot \partial_x^2 \bar{w}
\nonumber \\ 
&& - \ \frac{(m-1)^2}{m}
     \Big(  \frac{ \psi'(\bar{w})}{\psi(\bar{w})} \Big)^{'}
        \cdot \psi \cdot  (u+\varepsilon)^{q-m-2} u
        \cdot \partial_x^2 v \cdot \partial_x \bar{w}
\nonumber \\
&& - \ \frac{m-1}{m} \cdot 
              \psi(\bar{w}) \cdot 
             (u +\varepsilon)^{q-m-2} 
             \Big((q-2) u +\varepsilon
             \Big) \cdot  \partial_x^2 v \cdot \partial_x \bar{w}
\nonumber \\
&& - \ \frac{(m-1)^2}{m} \cdot  \frac{ \psi(\bar{w})^2}{\psi'(\bar{w})} 
        \cdot (u +\varepsilon)^{q-m-2} u
\cdot \partial_x^3 v.
\end{eqnarray}
Let us put $U:=|\partial_x \bar{w}|^2$,
then the following identities hold.
\begin{eqnarray} \label{der;k1..} 
\partial_x \bar{w} \cdot \partial_x^2 \bar{w} 
       = \frac12  \partial_x U,
\qquad
\partial_x \bar{w} \cdot \partial_x^3 \bar{w} 
     & = & 
       \frac{1}{2} \partial_x^2 U - (\partial_x^2 \bar{w})^2. 
\end{eqnarray}
Multiplying (\ref{bar u 1..}) by $\partial_x \bar{w}$
and using (\ref{der;k1..}), the resultant equation in terms of $U$ reads: 
\begin{eqnarray} \label{U..}
\frac{1}{2} \cdot \partial_t U
& = & \Big((m-1) \psi \big( \frac{\psi^{''}}{\psi^{'}} \big)^{'} 
         + m\psi^{''} 
      \Big) U^2
     + \Big( 
          (m-1) \cdot \psi \cdot \frac{\psi^{''}}{\psi^{'}} 
                                  + \frac{m+1}{2} \psi^{'} 
       \Big) \partial_x \bar{w} \cdot \partial_x U
\nonumber \\
&& 
+ \ (m-1) \psi 
             \cdot \Big( 
                 \frac12 \partial_x^2 U - (\partial_x^2 \bar{w})^2  
                   \Big)
\nonumber \\
&& - \ (q-2) \cdot \frac{1}{m} \cdot \psi^{'}(\bar{w})
         \cdot (u+\varepsilon)^{q-m-2} \cdot 
                \Big((q-1)u+2\varepsilon \Big) 
         \cdot \partial_x v   
         \cdot \partial_x \bar{w} \cdot U
\nonumber \\
&& - \ (q-2) \cdot \frac{m-1}{m} \cdot \psi(\bar{w})
         \cdot (u+\varepsilon)^{q-m-2} u 
         \cdot \partial_x^2 v 
         \cdot  U
\nonumber \\
&& - \ (q-2) \cdot \frac{m-1}{m} \cdot \psi(\bar{w})
         \cdot (u+\varepsilon)^{q-m-2}  u
         \cdot \partial_x v \cdot 
             \frac{1}{2} \partial_x U 
\nonumber 
\end{eqnarray}
\begin{eqnarray}
&& - \   \frac{m-1}{m} \cdot \psi(\bar{w})
         \cdot (u+\varepsilon)^{q-m-1}  
         \cdot \partial_x^2 v \cdot U 
\nonumber \\
&& - \   \frac{m-1}{m} \cdot \psi(\bar{w})
         \cdot (u+\varepsilon)^{q-m-1}  
         \cdot \partial_x v  \cdot \frac{1}{2} \partial_x U 
\nonumber \\ 
&& - \ \frac{(m-1)^2}{m}
     \Big(  \frac{ \psi'(\bar{w})}{\psi(\bar{w})} \Big)^{'}
        \cdot \psi \cdot  (u+\varepsilon)^{q-m-2} u
        \cdot \partial_x^2 v \cdot U
\nonumber \\
&& - \ \frac{m-1}{m} \cdot 
              \psi(\bar{w}) \cdot 
             (u +\varepsilon)^{q-m-2} 
             \Big((q-2) u +\varepsilon
             \Big) \cdot  \partial_x^2  v \cdot U
\nonumber \\
&& - \ \frac{(m-1)^2}{m} \cdot  \frac{ \psi(\bar{w})^2}{\psi'(\bar{w})} 
        \cdot (u +\varepsilon)^{q-m-2} u
\cdot \partial_x^3 v \cdot \partial_x \bar{w}.
\end{eqnarray}
\qquad
We consider a sequence $\{\eta_k(x)\}_{k=-\infty}^{\infty}$ 
of cut-off functions such that
\begin{eqnarray}
&& {\rm supp} \ \eta_k  =  \{x \in \R; \ -2+k \le x \le 2+k \}, 
\label{ino-oi}
\\
&& \eta_k(x)  = 1 \qquad {\rm for} \ -1+k \le x \le 1+k,
\label{ino-oi-222}
\end{eqnarray}
with
\begin{eqnarray}\label{opsi}
&& |\partial_x \eta_k(x)| 
 \le  c_1 (\eta_k(x))^{\frac{3}{4}},
\quad
-c_2 \eta_k(x)
\ \le \ 
\partial_x^2 \eta_k(x) 
\ \le \ c_3
\qquad {\rm for} \ x \in \R,
\end{eqnarray}
where $c_1, c_2$ and $c_3$ are positive constants independent of $k$.
In Remark 5 below, we give an example of such $\{\eta_k\}_{k=-\infty}^{\infty}$. 
\vspace{2mm} \\
\qquad
Multiplication of (\ref{U..}) by $\eta_k$ yields 
\begin{eqnarray} \label{U..ij}
\frac{1}{2} \cdot \partial_t (\eta_k U)
& = & 
\Big((m-1) \psi \big( \frac{\psi^{''}}{\psi^{'}} \big)^{'} 
         + m\psi^{''} 
      \Big) \eta_k U^2 
\ - \ (m-1) \psi 
             (\partial_x^2 \bar{w})^2  \eta_k
\nonumber \\
&& \ + \ I \cdot \partial_x^2(\eta_k U)
\ + \ J \cdot \partial_x(\eta_k U)
\ + \ R_k,
\end{eqnarray}
where
\begin{eqnarray*}
I & := &  \frac{m-1}{2} \psi, 
\nonumber \\
J & := & \Big( 
          (m-1) \cdot \psi \cdot \frac{\psi^{''}}{\psi^{'}} 
                                  + \frac{m+1}{2} \psi^{'} 
       \Big) \partial_x \bar{w}
\ - \ (m-1) \psi \cdot \partial_x \eta_k
\nonumber \\
&& - \frac{m-1}{2m} \cdot \psi (u+\varepsilon)^{q-m-2}\Big( (q-1)u+\varepsilon \Big)
              \partial_x v,
\end{eqnarray*}
and $R_k$ is regarded as the remainder term defined by 
\begin{eqnarray} \label{RK}
R_k & := & 
\sum_{j=1}^7 R_k^{(j)}
\end{eqnarray}
with
\begin{eqnarray*}
R_k^{(1)} & := & - \Big( 
          (m-1) \cdot \psi \cdot \frac{\psi^{''}}{\psi^{'}} 
                                  + \frac{m+1}{2} \psi^{'} 
       \Big) \partial_x \eta_k \cdot \partial_x \bar{w} \cdot U
\nonumber \\
&& 
- \ \frac{q-2}{m} \cdot \psi^{'}(\bar{w})
         \cdot (u+\varepsilon)^{q-m-2} \cdot 
                \Big((q-1)u+2\varepsilon \Big) 
         \cdot \partial_x v \cdot \eta_k
         \cdot \partial_x \bar{w} \cdot U,
\nonumber \\
R_k^{(2)} & := & 
(m-1) \cdot \psi \cdot (\partial_x \eta_k)^2 U,
\nonumber \\
R_k^{(3)} & := & \frac{m-1}{2m} \cdot \psi(\bar{w})
         \cdot (u+\varepsilon)^{q-m-2}  \Big( (q-1)u + \varepsilon \Big)
         \cdot \partial_x v \cdot \partial_x \eta_k \cdot U, 
\nonumber \\
R_k^{(4)} & := & 
    - \ \frac{m-1}{m} \cdot \psi(\bar{w})
         \cdot (u+\varepsilon)^{q-m-2} \Big( (2q-3)u + 2\varepsilon \Big) 
         \cdot \partial_x^2 v \cdot \eta_k 
           U,
\nonumber \\
R_k^{(5)} & := & 
- \ \frac{m-1}{2} \cdot \psi \cdot \partial_x^2 \eta_k \cdot U,
\nonumber \\
R_k^{(6)} & := & 
- \ \frac{(m-1)^2}{m} \cdot
     \Big(  \frac{ \psi'(\bar{w})}{\psi(\bar{w})} \Big)^{'}
        \cdot \psi \cdot  (u+\varepsilon)^{q-m-2} u
        \cdot \partial_x^2 v \cdot \eta_k U,
\nonumber \\
R_k^{(7)} & := & 
- \ \frac{(m-1)^2}{m} \cdot  \frac{ \psi(\bar{w})^2}{\psi'(\bar{w})} 
        \cdot (u +\varepsilon)^{q-m-2} u
\cdot \partial_x^3 v \cdot \partial_x \bar{w} \cdot \eta_k.
\end{eqnarray*}
Now we choose the transformation $\psi(r)$ by
\begin{eqnarray}\label{PSIAT}
\psi(r) := \frac{m}{m-1} (\|u_0\|_{L^{\infty}(\R)}+2+\varepsilon)^{m-1} 
                  \cdot \frac{r}{3}(4-r),
\qquad 0 \le r \le 1.
\end{eqnarray}
Then we observe that 
the coefficient of the first term $\eta_k U^2$ of the right-hand side in
(\ref{U..ij}) is negative, in particular
\begin{eqnarray} \label{psi;negative}
\Big((m-1) \psi \big( \frac{\psi^{''}}{\psi^{'}} \big)^{'} 
         + m\psi^{''} 
      \Big) \eta_k U^2
& \le &  - M \cdot \eta_k U^2
\qquad  
{\rm for} \ (x,t) \in \R \times [0,T_0),
\end{eqnarray}
where 
$ \displaystyle M:=\frac{m(11m - 3)}{12(m-1)} >0$.
\vspace{2mm} \\
Indeed, since Proposition \ref{pr;apto} states that
$$
\frac{m}{m-1} \varepsilon^{m-1}
 \le w(x,t) 
             \le \frac{m}{m-1}
(\|u_0\|_{L^{\infty}(\R)}+2+\varepsilon)^{m-1} \ =: \ L
$$ 
holds for all $(x,t) \in \R \times [0,T_0)$, 
the definition $\psi(\bar{w})=w$ yields
\begin{eqnarray} \label{re;barw}
(0<) \ 2-\sqrt{4-\frac{3m\varepsilon^{m-1}}{(m-1)L}} 
\le \bar{w}(x,t) \le 1,
\qquad (x,t) \in \R \times [0,T_0)
\end{eqnarray}
for sufficiently small $\varepsilon>0.$
Moreover, by (\ref{re;barw}) we have
\begin{eqnarray}
&& 
\frac{2L}{3} 
\ \le \ 
\psi^{'}(\bar{w}) 
\ = \ \frac{2L}{3}(2-\bar{w})
\ \le \ \frac{4L}{3},
\quad
\psi^{''}(\bar{w}) 
\ = \ -\frac{2L}{3},
\label{w1} \\
&& 
\frac{1}{2}
\ \le \ 
\left|
\frac{\psi^{''}}{\psi^{'}}
\right|
\ \le \ 1,
\quad
-1 \ \le \ 
\Big( \frac{\psi^{''}}{\psi^{'}}  \Big)^{'}
\ = \ \frac{\psi^{'}\psi^{'''}- (\psi^{''})^2}{(\psi^{'})^2}
\ = \ -\Big( \frac{\psi^{''}}{\psi^{'}}   \Big)^2
\ \le \  - \frac14.
\label{w2}
\end{eqnarray}
Now from (\ref{w1}) and (\ref{w2}), 
we see that the left-hand side of (\ref{psi;negative}) is bounded by
\begin{eqnarray*}
\Big((m-1) \psi \big( \frac{\psi^{''}}{\psi^{'}} \big)^{'} 
         + m\psi^{''} 
      \Big) \cdot \eta_k U^2
& \le & 
\Big(- \frac{(m-1) L}{4} - \frac{2mL}{3}
      \Big) \cdot \eta_k U^2
\nonumber \\
& = &  - \frac{m(11m - 3)}{12(m-1)} \cdot \eta_k U^2 
= -M \cdot \eta_k U^2 <0.
\end{eqnarray*}
On the other hand, 
suppose that $\eta_k U$ attains its maximum at the point 
$(x_0,t_0) \in \R \times (0,T_0)$.
Then it holds by $\eta_k \ge 0$ that
\begin{eqnarray} \label{ux,ut}
\frac{1}{2} \cdot \partial_t (\eta_k U)(x_0,t_0) \ \ge \ 0,
\quad
\partial_x^2 (\eta_k U)(x_0,t_0) \ \le \ 0,
\quad \partial_x (\eta_k U)(x_0,t_0) \ = \ 0.
\end{eqnarray}
Combining (\ref{psi;negative}), (\ref{w1})--(\ref{ux,ut}) 
                  with (\ref{U..ij}), we obtain 
\begin{eqnarray}\label{uiui}
M \cdot \eta_k U^2
& \le & \sum_{j=1}^7 R_k^{(j)},
\label{IJ1-J5}
\end{eqnarray}
where $\{ R_k^{(j)} \}_{1 \le j \le 7}$ is given by (\ref{RK}).
\vspace{2mm} \\
\qquad
We are going to estimate the seven terms $\{R_k^{(j)}\}_{1 \le j \le 7}$.
To this end, firstly  
integrating the second equation of (KS)$_{\varepsilon}$ on $(-\infty, x)$,
we have
\begin{eqnarray*}
\partial_x v(x,t) 
     = \gamma \int_{-\infty}^x v(y,t) \ dy 
               - \int_{-\infty}^x u(y,t) \ dy, 
\end{eqnarray*}
which yields
\begin{eqnarray}\label{proof;th-2}
\sup_{0<t<T_0} \| \partial_x v(t)\|_{L^{\infty}(\R)} 
& \le & 2\| u_0 \|_{L^1(\R)}.
\end{eqnarray}
Here we have used the fact that 
$v(x,t) > 0$ together with 
$
\gamma \int_{-\infty}^{\infty} v(y,t) \ dy= 
        \int_{-\infty}^{\infty} u(y,t) \ dy$ for all $t \in [0,T_0)$.
\vspace{2mm} \\
\qquad
By (\ref{opsi}), (\ref{proof;th-2}), Young's inequality, 
         and the relation $\partial_x^2 v = \gamma v-u$, 
we have for $q \ge m+1$ that
\begin{eqnarray} \label{J1-J4}
\sum_{j=1}^4 R_k^{(j)} & \le & 
C \ + \ \frac{M}{8} \cdot \eta_k U^2,
\end{eqnarray}
where
$C$ is a constant depending on $m,\gamma,q$ and $u_0$.
By (\ref{opsi}) and Young's inequality, we have
\begin{eqnarray}\label{sumgni}
R_k^{(5)}
& \le & 
\frac{c_2(m-1)}{2} \psi \cdot \eta_k U
\ \le \ 
C \ + \ \frac{M}{8} \cdot \eta_k U^2.
\end{eqnarray}
\qquad
We are now going to estimate $R_k^{(6)}$.
Since
\begin{eqnarray*}
-\frac{2L}{3} \cdot \frac{1}{\psi} 
-\Big( \frac{4L}{3}  \Big)^2
\cdot  \frac{1}{\psi^2} 
\ \le \ 
\Big( \frac{\psi^{'}}{\psi}  \Big)^{'}
\ = \ \frac{\psi\psi^{''}- (\psi^{'})^2}{(\psi)^2}
\ \le \ -\frac{2L}{3} \cdot \frac{1}{\psi}
\ < \ 0,
\label{w3}
\end{eqnarray*}
and since $\partial_x^2 v = \gamma v-u$, 
by the hypothesis that $q \ge 2m$, we have that
\begin{eqnarray} 
R_k^{(6)} & = & - \frac{(m-1)^2}{m} \cdot 
            \Big(\frac{ \psi'(\bar{w})}{\psi(\bar{w})} \Big)^{'} u      
\cdot 
\psi(\bar{w}) 
         \cdot (u+\varepsilon)^{q-m-2} 
         \cdot \partial_x^2 v 
         \cdot \eta_k U
\nonumber \\
& \le & 
\frac{(m-1)^2}{m} \cdot 
\frac{2L}{3} 
         \cdot (u+\varepsilon)^{q-m-1} 
         \cdot \gamma v 
         \cdot  \eta_k U
\nonumber \\
&& \ + \ (m-1) \cdot 
\Big(\frac{4L}{3} \Big)^2  
         \cdot (u+\varepsilon)^{q-2m} 
         \cdot \gamma v 
         \cdot  \eta_k U
\nonumber \\
& \le & 
C \ + \ \frac{M}{8} \cdot \eta_k U^2.
\label{JJ}
\end{eqnarray}
Indeed, since 
$v(x,t)$ satisfies
\begin{eqnarray} \label{semi; z}
v(x,t)
& = & \int_{\R} G(x-y) \cdot u(y,t) \ dy 
\end{eqnarray} 
with
the Bessel potential $G(x)$ 
                  which can be express as
\begin{eqnarray} \label{GI}
G(x) & = & \frac{1}{\sqrt{4\pi}} 
       \int_0^{\infty} s^{-\frac{1}{2}} \cdot e^{-\gamma s-\frac{|x|^2}{4s}} \ ds
\ = \ \frac{e^{-\sqrt{\gamma}|x|}}{2 \gamma},
\end{eqnarray}
it holds that
$
G \in L^p(\R)
$ for all $1 \le p \le \infty$ 
and that
\begin{eqnarray*}
\sup_{0<t<T} \| v(t) \|_{L^{\infty}(\R)} 
& \le & C \sup_{0<t<T} \|u(t)\|_{L^1(\R)} \ = \ C \|u_{0\varepsilon}\|_{L^1(\R)} 
\ \le \ \|u_0\|_{L^1(\R)}.
\end{eqnarray*}
By (\ref{psi;w-1}) and (\ref{psi;w-2}), it holds
\begin{eqnarray} \label{J5}
\lefteqn{
\psi \cdot (u+\varepsilon)^{q-m-2} u \cdot \partial_x^3 v
} \nonumber \\
& = & 
\frac{m}{m-1} (u+\varepsilon)^{q-3} u \cdot 
                           (\gamma \partial_x v - \partial_x u )
\nonumber \\
& \le & 
      \frac{m}{m-1} \cdot \gamma \cdot 
                               (u+\varepsilon)^{q-2} \cdot 2\| u_0 \|_{L^1}  
        +
      \frac{1}{m-1} \cdot \psi^{'}(\bar{w}) (u+\varepsilon)^{q-m-1} u\cdot 
                        \partial_x \bar{w},
\end{eqnarray}
which yields
\begin{eqnarray} \label{J5-3}
R_k^{(7)} & \le & 
C \ + \ \frac{M}{8} \cdot \eta_k U^2.
\end{eqnarray} 
Substituting (\ref{J1-J4}), (\ref{sumgni}), 
                (\ref{JJ}) and (\ref{J5-3}) into (\ref{IJ1-J5}),
we obtain
\begin{eqnarray*}
M \cdot \eta_k U^2
& \le & C.
\end{eqnarray*} 
Recalling $U=|\partial_x \bar{w}|^2$, 
we have by (\ref{ino-oi}) and the above estimate that
\begin{eqnarray} \label{Rm}
|\partial_x \bar{w}|^2 \ =: \ U
& \le &
C \ \qquad {\rm for} \ -1+k \le x \le 1+k, \quad 0 \le t <T_0, 
\end{eqnarray}
where $C$ is a constant independent of $\varepsilon$ and $k$.
Repeating the same argument as the above for $k=0, \pm 1, \pm 2, \cdots,$
we obtain the upper bound of $|\partial_x \bar{w}|$ 
which is independent of $\varepsilon$ in the whole interval $\R$.
\vspace{2mm} \\
\qquad
We recall the definition of $w$ and $\psi({\bar{w}})$:
\begin{eqnarray}\label{relation;u-w-barw}
\frac{m}{m-1}(u_{\varepsilon}+\varepsilon)^{m-1}
& = & 
w 
\ = \ 
\psi(\bar{w}).
\end{eqnarray}
Differentiating both sides of (\ref{relation;u-w-barw}) with respect to $x$,
we have by (\ref{w1}) and (\ref{Rm}) that
\begin{eqnarray*}\label{fga}
\left| \partial_x (u_{\varepsilon}+\varepsilon)^{m-1} \right|
& = & 
\frac{m-1}{m}
\cdot 
\psi^{'}(\bar{w}) \cdot 
|\partial_x \bar{w}|
\ \le \  C \quad {\rm for \ all} \ (x,t) \in \R \times [0,T_0),
\end{eqnarray*}
which yields (\ref{adj}). This proves Lemma \ref{lip}.
\vspace{2mm} \\
{\bf Remark 5.}
In the proof of Lemma \ref{lip}, we have used a sequence
$\{\eta_k(x)\}_{k=-\infty}^{\infty}$ of cut-off functions with properties 
(\ref{ino-oi})--(\ref{opsi}).
Taking $\eta(x)$ as
\begin{eqnarray*}
&& 
\eta(x) = \left\{
\begin{array}{llll}
 0 & \quad {\rm for} \quad x \le -2 \\
                       2(2+x)^4 & 
                         \quad {\rm for} \quad -2 < x < -\frac{3}{2}, \\
                       1-2(x+1)^4 & 
                         \quad {\rm for} \quad -\frac{3}{2} < x \le -1, \\
                    1 & \quad {\rm for} \quad -1 \le x \le 1, \\
                       1-2(x-1)^4 & 
                         \quad {\rm for} \quad  1 < x \le \frac{3}{2}, \\
                       2(2-x)^4 & 
                         \quad {\rm for} \quad \frac{3}{2} < x < 2, \\
                    0 & \quad {\rm for} \quad x \ge 2, 
\end{array}
                 \right.
\end{eqnarray*}
and then defining $\eta_{k}$ by $\eta_{k}(x) := \eta(x-k)$ 
for $k=0, \pm 1, \pm 2, \cdots$,
we see that $\{\eta_k(x)\}_{k=-\infty}^{\infty}$ has the desired
properties (\ref{ino-oi})--(\ref{opsi}).
\ \\
\section{Proof of Theorem \ref{thm:lip}}
\label{section;3}
\qquad
Let us first show that for every $t \in [0,T_0)$,
$\{u_{\varepsilon}(\cdot,t)\}_{\varepsilon>0}$ is a sequence of
uniformly bounded and equi-continuous functions in $\R$.
Indeed, the uniform bound is a consequence of (\ref{estK0}).
By (\ref{estK0}), (\ref{adj}) and (\ref{holder-est}) with $u$ replaced
by $u_{\varepsilon}+\varepsilon$,
it holds
\begin{eqnarray*}
|u_{\varepsilon}(x,t)-u_{\varepsilon}(y,t)| 
& \le & C(\|u_0\|_{L^{\infty}}+2)^2 |x-y|^{\mu},
\quad
\mu=\min\{1,\frac{1}{m-1}\}
\end{eqnarray*} 
for all $x,y \in \R, \ 0 \le t <T_0$, and all $\varepsilon>0$,
where $C$ is the same constant as in (\ref{adj}).
This implies that $\{u_{\varepsilon}(\cdot,t)\}_{\varepsilon>0}$ is a
family of equi-continuous functions in $\R$ for all $0 \le t <T_0$.
Hence by the Ascoli-Arzela theorem, 
there is a subsequence of $\{u_{\varepsilon}(\cdot,t)\}_{\varepsilon>0}$,
which we denoted by $\{u_{\varepsilon}(\cdot,t)\}_{\varepsilon>0}$
itself such that
\begin{eqnarray}\label{Lm;strong--sta}
u_{\varepsilon}(\cdot,t) & \longrightarrow & u(\cdot,t) 
\qquad {\rm as} \ \ \varepsilon \to 0
\end{eqnarray}
uniformly in every compact interval $I \subset \R$.
\vspace{2mm} \\
\qquad
On the other hand, by (\ref{adj}) and
the weakly-star compactness of $L^{\infty}(Q_{T_0})$,
there exists a sequence of $\{u_{\varepsilon}\}_{\varepsilon>0}$,
which we denote by $\{u_{\varepsilon}\}_{\varepsilon>0}$ itself for simplicity,
and a function $\tilde{u} \in L^{\infty}(Q_{T_0})$ such that
\begin{eqnarray*}
\partial_x (u_{\varepsilon}+\varepsilon)^{m-1} \ 
\rightarrow 
\tilde{u} \hspace{1cm} 
{\rm weakly-star}  
&              \ {\rm in} \ 
L^{\infty}(Q_{T_0}) 
\end{eqnarray*}
with
\begin{eqnarray*}
\| \tilde{u} \|_{L^{\infty}(Q_{T_0})} & \le & 
\liminf_{\varepsilon \to +0} 
\| \partial_x (u_{\varepsilon}+\varepsilon)^{m-1} \|_{L^{\infty}(Q_{T_0})}.
\end{eqnarray*}
By (\ref{Lm;strong--sta}), it is easy to see that $\tilde{u}= \partial_x u^{m-1}$,
which yields the desired estimate (\ref{pr-adj}).
\vspace{2mm} \\
\qquad 
Next, we shall show that
$\partial_x u^{m-1+\delta}(\cdot,t)$ is a continuous function in $\R$ 
for all $0<t<T_0$ and for all $\delta>0$ with the additional property that
$\partial_x u^{m-1+\delta}(x,t)=0$ at the point $(x,t)$ such as $u(x,t)=0$.
To this aim, we follow a similar argument employed in Aronson \cite{Ar3}.
Let $u(x_0,t_0)>0$. 
Then we see by the standard argument that 
both $\partial_x u$ and $\partial_x u^{m-1+\delta}$ with $\delta>0$ 
    are continuous functions in a neighbourhood of $(x_0,t_0)$. 
Therefore, it suffices to prove that 
$\partial_x u^{m-1+\delta}(\cdot,t)$ is a continuous function in a
neighbourhood of $x_1$ such as $u(x_1,t)=0$ 
with the additional property that $\partial_x u^{m-1+\delta}(x_1,t)=0$.
By virtue of (\ref{Lm;strong--sta}), 
for every $t \in [0,T_0)$ and every compact interval $I \subset \R$,
it holds that $u_{\varepsilon}(\cdot,t) \to u(\cdot,t)$ uniformly on $I$.             
Therefore, by Remark 2, there exists $a_0>0$ such that
\begin{eqnarray} \label{hani}
0 \le u_{\varepsilon}(x,t) 
& \le & 
|u_{\varepsilon}(x,t) -u(x,t)| + |u(x,t)-u(x_1,t)| +u(x_1,t)
\nonumber \\
& \le &  2 a^\mu 
\end{eqnarray}
holds for all $x \in I_a(x_1):=\{x \in \R; |x-x_1| < a\}$ 
and for all $0<a \le a_0$ and for all $0<\varepsilon<1$,
where $\mu:= \min \{1, \frac{1}{m-1} \}$.
\vspace{2mm} \\
\qquad
On the other hand, since we have
\begin{eqnarray}\label{sekibun}
u_{\varepsilon}^{m-1+\delta}(x,t)-u_{\varepsilon}^{m-1+\delta}(x^{\prime},t) 
& = & 
\frac{m-1+\delta}{m-1} \int_{x^{\prime}}^x u_{\varepsilon}^{\delta}(x,t)
                              \cdot 
                            \partial_x u_{\varepsilon}^{m-1}(x,t) \ dx,
\end{eqnarray}
it follows from (\ref{hani}),(\ref{sekibun}) and Lemma \ref{lip} that
\begin{eqnarray}\label{ilet}
|u_{\varepsilon}^{m-1+\delta}(x,t)-u_{\varepsilon}^{m-1+\delta}(x^{\prime},t)| 
& \le & 
C (2a^{\mu})^{\delta} |x-x^{\prime}|
\qquad
{\rm for \ all} \ x,x^{\prime} \in I_a(x_1)
\end{eqnarray} 
and for all  $0 < a \le a_0$ and for all $0<\varepsilon<1$,
where $C$ depends on $m,\gamma,q,u_0$ but not on $\varepsilon$.
Letting $\varepsilon \to +0$ in (\ref{ilet}), 
    we have by (\ref{Lm;strong--sta}) that
\begin{eqnarray} \label{base}
\lefteqn{
|u^{m-1+\delta}(x,t)-u^{m-1+\delta}(x^{\prime},t)| 
} \nonumber \\
& \le & 
C (2a^{\mu})^{\delta} |x-x^{\prime}|
\quad
{\rm for \ all} \ x,x^{\prime} \in I_a(x_1)
\ {\rm and \ all} \ 0 < a \le a_0.
\end{eqnarray} 
Taking $x=x_1$ in (\ref{base}) and then letting $x^{\prime} \to x_1$, we have
\begin{eqnarray*}
|\partial_x u^{m-1+\delta}(x_1,t)| & \le & C(2a^{\mu})^{\delta}, \quad 0<a \le a_0.
\end{eqnarray*}
Hence we have by letting $a \to 0$ that
\begin{eqnarray*}
\partial_x u^{m-1+\delta}(x_1,t)=0.
\end{eqnarray*}
Similarly, letting $x^{\prime} \to x$ in (\ref{base}), we have
\begin{eqnarray*}\label{base-2}
|\partial_x u^{m-1+\delta}(x,t)| & \le & C(2a^{\mu})^{\delta}
\qquad {\rm for \ all} \ 0<a\le a_0,
\end{eqnarray*}
which implies that $\partial_x u^{m-1+\delta}(\cdot,t)$ is continuous at $x_1$.
Since $x_1$ can be taken arbitrary in such a way that $u(x_1,t)=0$,
we conclude that 
$\partial_x u^{m-1+\delta}(\cdot,t)$ 
is a continuous function in $\R$ for all $t \in [0,T_0)$
with the additional property that
$\partial_x u^{m-1+\delta}(x,t)=0$ for the point $(x,t)$ such as $u(x,t)=0$.
\vspace{2mm} \\
\qquad
The case of $1<m<2$ can be handled in a similar manner as above 
and we conclude that 
$\partial_x u(\cdot,t)$ is a continuous function in $\R$ for all $t \in [0,T_0)$
with the additional property that
$\partial_x u(x,t)=0$ for the point $(x,t)$ such as $u(x,t)=0$.
This completes the proof of Theorem \ref{thm:lip}.
\section{Proof of Theorem \ref{thm;finite speed} }
\label{section;4}
\qquad
Let $(u_{\varepsilon},v_{\varepsilon})$ 
be the unique strong solution of (KS)$_{\varepsilon}$
given by Proposition \ref{pr;apto}.
For a fixed $R>0$, we take $a,b>0$ such as $-R<a<b<R$ 
and consider the following ordinary differential equations:
\begin{equation*}
{\rm (IE)_{\xi}} :\left\{
    \begin{array}{llll}
         & \xi_{\varepsilon}^{'}(t)  = 
\frac{m}{m-1} \partial_x \big( u_{\varepsilon} + \varepsilon
                     \big)^{m-1} (\xi_{\varepsilon}(t),t)
 - \big(u_{\varepsilon} 
              + \varepsilon
                  \big)^{q-3} 
               u_{\varepsilon}
              \cdot \partial_x v_{\varepsilon}(\xi_{\varepsilon}(t),t), 
\ \ 0 \le t<T_0,
\\
         & \xi_{\varepsilon}(0)   = a,
     \end{array}
  \right.
\end{equation*}
and
\begin{equation*}
{\rm (IE)_{\Xi}} :\left\{
    \begin{array}{llll}
         & \Xi_{\varepsilon}^{'}(t)  = 
\frac{m}{m-1} \partial_x \big( u_{\varepsilon} + \varepsilon
                     \big)^{m-1} (\Xi_{\varepsilon}(t),t)
 - \big(u_{\varepsilon} 
              + \varepsilon
                  \big)^{q-3} 
               u_{\varepsilon}
              \cdot \partial_x v_{\varepsilon}(\Xi_{\varepsilon}(t),t), 
\ \ 0 \le t<T_0,
\\
         & \Xi_{\varepsilon}(0)   = b.
     \end{array}
  \right.
\end{equation*}
By Remark 4 (iii), (\ref{estK0}), (\ref{adj}) and (\ref{proof;th-2}), 
     we have
\begin{eqnarray}\label{2R}
\partial_x (u_{\varepsilon}+\varepsilon)^{m-1}
\in C^1([-2R,2R] \times [0,T_0))
\end{eqnarray}
and
\begin{eqnarray*}
&& \sup_{0<\varepsilon<1}
\Big(
\sup_{0<t<T_0}
\{   
  \| \partial_x 
        (u_{\varepsilon}+\varepsilon)^{m-1}(\cdot,t) \|_{L^{\infty}(-2R,2R)}
+
  \|
\big(u_{\varepsilon} 
              + \varepsilon
                  \big)^{q-3} 
               u_{\varepsilon}
              \cdot \partial_x v_{\varepsilon}(\cdot,t) 
 \|_{L^{\infty}(-2R,2R)}
\}
\Big)
\nonumber \\
&& \quad \le  C,
\end{eqnarray*}
where $C=C(m,\gamma,q,u_0)$.
We now chose $R>0$ large enough such that $\frac{2R}{C} > T_0$.
Then, it follows from the well-known theorem on the existence and 
uniqueness of local solutions to the initial value problem for the ordinary
differential equations that
both (IE)$_{\xi}$ and (IE)$_{\Xi}$ have 
a unique $C^1$-solution $\xi_{\varepsilon}(t)$ and $\Xi_{\varepsilon}(t)$
  on $[0,T_0)$ for all $\varepsilon>0$, respectively.
\vspace{2mm} \\
We consider the following domain: 
\begin{eqnarray*}
D_{\tau} 
& := & 
\displaystyle 
\bigcup_{t \in [0,\tau]} 
I_{t} \times \{t\}, \quad 
I_{t} \ := \  
\Big\{ x \in \R; \ 
        \xi_{\varepsilon}(t) \le x \le \Xi_{\varepsilon}(t)
          \Big\}
\quad \ {\rm for} \ 0 <\tau < T_0.
\end{eqnarray*}
By the local uniqueness 
  of the initial value problem (IE)$_\xi$ and (IE)$_{\Xi}$, 
    we obtain that
$$
\xi_{\varepsilon}(t) < \Xi_{\varepsilon}(t)
\qquad
{\rm for \ all } \
0 \le t <T_0.
$$
Let us define the gradient $\overrightarrow{\nabla}$ and 
the vector $\textbf{F}$ on $(x,t)$ by
\begin{eqnarray*}
\overrightarrow{\nabla} 
& := & (\partial_x, \partial_t),
\quad
\textbf{F}(x,t) 
\ := \ 
\Big(- \partial_x  (u_{\varepsilon}+\varepsilon)^m 
           + (u_{\varepsilon} + \varepsilon)^{q-2} u_{\varepsilon}                              
                    \cdot \partial_x v_{\varepsilon},
  \ \ u_{\varepsilon}+\varepsilon
\Big).
\end{eqnarray*}
Then it follows from the first equation of (KS)$_{\varepsilon}$ that
\begin{eqnarray} \label{a1}
\lefteqn{
\int_{D_{\tau}} \overrightarrow{\nabla} \cdot \textbf{F}(x,t) \ dxdt
} \nonumber \\
& = & \int_{D_{\tau}} \partial_t u_{\varepsilon} 
    - \partial_x \Big( 
        \partial_x (u_{\varepsilon} +\varepsilon)^m 
           - (u_{\varepsilon}  + \varepsilon)^{q-2} u_{\varepsilon}                               
                    \cdot   \partial_x v_{\varepsilon}  
             \Big) \ dxdt 
\ = \ 0
\end{eqnarray}
for all $0<\tau<T_0$.
Taking two curves $C_1$ and $C_2$ as
\begin{eqnarray*}
C_1 := \{(x,t) = (\xi_{\varepsilon}(t),t) ; \ 0<t<\tau \},
\quad
C_2 := \{(x,t) = (\Xi_{\varepsilon}(t),t) ; \ 0<t<\tau \},
\end{eqnarray*}
we have
\begin{eqnarray*}
\partial D_{\tau} & = & I_0 \ \cup \ C_1 \ \cup \ C_2 \ \cup \ I_{\tau}.
\end{eqnarray*}
Hence, the Stokes formula gives
\begin{eqnarray} \label{a2}
\lefteqn{
0= \int_{D_{\tau}} \overrightarrow{\nabla} \cdot \textbf{F}(x,t) \ dxdt
} \nonumber \\
& = & 
\int_{\partial D_{\tau}} \textbf{F}(x,t) \cdot {\bf n} \ dS
\nonumber \\
& = & \int_a^b \textbf{F}(x,0) \cdot (0,-1) \ dx
           - \int_{\Xi_{\varepsilon}(\tau)}^{\xi_{\varepsilon}(\tau)} 
               \textbf{F}(x,t) \cdot (0,1) \ dx
           + \int_{C_1} 
                         \textbf{F} \cdot {\bf n_1} \ dS
           + \int_{C_2} 
                         \textbf{F} \cdot {\bf n_2} \ dS
\nonumber \\
\quad
& = & - \int_a^b (u_{0\varepsilon}+\varepsilon) \ dx
           + \int_{\xi_{\varepsilon}(\tau)}^{\Xi_{\varepsilon}(\tau)}
                     (u_{\varepsilon}+\varepsilon)  \ dx
           + \int_{C_1} 
                         \textbf{F} \cdot {\bf n_1} \ dS
           + \int_{C_2} 
                         \textbf{F} \cdot {\bf n_2} \ dS,
\end{eqnarray}
where ${\bf n_1}$ and ${\bf n_2}$ denote the unit outer normals to $C_1$
and $C_2$, respectively.
Since
\begin{eqnarray*}  
{\bf n_1} & = & \frac{(1, \xi_{\varepsilon}^{\prime}(t))}
                     {\sqrt{1+(\xi_{\varepsilon}^{\prime}(t))^2}},
\qquad
{\bf n_2} \ = \ \frac{(1, \Xi_{\varepsilon}^{\prime}(t))}
                     {\sqrt{1+(\Xi_{\varepsilon}^{\prime}(t))^2}},
\end{eqnarray*}
we have by (IE)$_{\xi}$ and (IE)$_{\Xi}$ that
\begin{eqnarray}\label{agh-1}
\textbf{F} \cdot {\bf n_1} & = & 0
\quad {\rm on} \ \ C_1,
\qquad
\textbf{F} \cdot {\bf n_2} \ = \ 0
\quad {\rm on} \ \ C_2.
\end{eqnarray}
Combining (\ref{a1})--(\ref{agh-1}), we have
\begin{eqnarray} \label{key;ineq}
\int_{\xi_{\varepsilon}(\tau)}^{\Xi_{\varepsilon}(\tau)} 
            (u_{\varepsilon}(x,\tau)+\varepsilon)  \ dx
& = &  
\int_a^b (u_{0\varepsilon}(x)+\varepsilon) \ dx,
\qquad 0 \le \tau <T_0.
\end{eqnarray}
On the other hand,
we obtain from (\ref{proof;th-2}), Proposition \ref{pr;apto} and Lemma \ref{lip} that 
\begin{eqnarray*}
\sup_{0<\varepsilon<1} \| \xi_{\varepsilon}\|_{L^{\infty}(0,T_0)}
& \le & a+
\Big(   
\frac{m}{m-1} \cdot C 
+ 2(\|u_0\|_{L^{\infty}}+2)^{q-2} \cdot \| u_0 \|_{L^1}
\Big) \cdot T_0,
\label{111}
\\
\sup_{0<\varepsilon<1} \|\xi_{\varepsilon}^{'}\|_{L^{\infty}(0,T_0)}
& \le &
\frac{m}{m-1} \cdot C 
+ 2(\|u_0\|_{L^{\infty}}+2)^{q-2} \cdot \| u_0 \|_{L^1}.
\label{222}
\end{eqnarray*}
Hence it follows by the Ascoli-Arzela theorem that
there exists a subsequence of $\{\xi_{\varepsilon}(t)\}$, 
still denoted by $\{\xi_{\varepsilon}(t)\}_{\varepsilon>0}$, 
and a function $\xi \in C^{0,1}[0,T_0)$ such that
\begin{eqnarray}\label{xi}
\xi_{\varepsilon}(t) \ \to \ \xi(t) \qquad 
                          {\rm as} \ \varepsilon \to 0
            \quad {\rm uniformly \quad for \ every} \ \ t \in [0,T_0).
\end{eqnarray}
Obviously, a similar argument to $\Xi_{\varepsilon}(t)$ also holds, and
there exist a subsequence of $\{\Xi_{\varepsilon}(t)\}_{\varepsilon>0}$, 
still denoted by $\{\Xi_{\varepsilon}(t)\}$, and $\Xi \in C^{0,1}[0,T_0)$ such that
\begin{eqnarray}\label{Xi}
\Xi_{\varepsilon}(t) \ \to \ \Xi(t) \qquad 
                          {\rm as} \ \varepsilon \to 0
            \quad {\rm uniformly \quad for \ every} \ \ t \in [0,T_0).
\end{eqnarray}
Since $u_0 \equiv 0$ on $[a,b]$, by letting $\varepsilon \to 0$ in
(\ref{key;ineq}), we have
\begin{eqnarray} \label{conv;ineq}
\int_{\xi(t)}^{\Xi(t)} u(x,t) \ dx
\ = \ 
\int_a^b u_0(x) \ dx \ = \ 0
\end{eqnarray}
for all $0 \le t <T_0$.
Indeed, we may assume 
\begin{eqnarray*}
-2R < \xi_{\varepsilon}(t) < \Xi_{\varepsilon}(t) < 2R
\qquad
{\rm for \ all } \ \varepsilon>0, \ {\rm and \ all} \ 0 \le t <T_0,
\end{eqnarray*}
where $R>0$ is the same as in (\ref{2R}).
Hence it follows from (\ref{Lm;strong--sta}), 
        (\ref{xi}), (\ref{Xi}) and Proposition \ref{pr;local existence} that
\begin{eqnarray*} \label{conv;ineq-2}
\lefteqn{
\Big|
\int_{\xi_{\varepsilon}(t)}^{\Xi_{\varepsilon}(t)} (u_{\varepsilon}+\varepsilon)  \ dx 
-
\int_{\xi(t)}^{\Xi(t)} u  \ dx 
\Big|
} \nonumber \\
& \le & 
\Big| \int_{\xi_{\varepsilon}(t)}^{\Xi_{\varepsilon}(t)}
                              (u_{\varepsilon} - u) \ dx \Big|
\ + \ 
\varepsilon \int_{\xi_{\varepsilon}(t)}^{\Xi_{\varepsilon}(t)}  \ dx
\ + \ 
\Big| \int_{\xi_{\varepsilon}(t)}^{\Xi_{\varepsilon}(t)}
                              u \ dx 
- \int_{\xi(t)}^{\Xi(t)}
                              u \ dx 
\Big|
\nonumber \\
& \le & 
\Big(\sup_{-2R \le x \le 2R} 
  |u_{\varepsilon}(x,t)-u(x,t)| + \varepsilon  \Big) 
  (\Xi_{\varepsilon}(t)-\xi_{\varepsilon}(t))
\nonumber \\
&&
\quad \ + \
\|u\|_{L^{\infty}(Q_{T_0})}
\Big( |\Xi_{\varepsilon}(t)-\Xi(t)|
     +
      |\xi_{\varepsilon}(t)-\xi(t)|
\Big) 
\nonumber \\
& \le & 
4R \Big(\sup_{-2R \le x \le 2R} 
  |u_{\varepsilon}(x,t)-u(x,t)| + \varepsilon  \Big) 
\ + \
\Big( 
\| u_0 \|_{L^{\infty}} +2 
\Big)
\Big( |\Xi_{\varepsilon}(t)-\Xi(t)|
     +
      |\xi_{\varepsilon}(t)-\xi(t)|
\Big)
\nonumber \\
& \to & 0 \qquad {\rm as} \ \varepsilon \to +0,
\end{eqnarray*}
which yields (\ref{conv;ineq}).
Since $u$ is non-negative in $\R \times [0,T_0)$,
we conclude from (\ref{conv;ineq}) that
\begin{eqnarray*}
u(x,t)=0 \qquad {\rm for} \quad \xi(t) \le x \le \Xi(t), \ 0 \le t<T_0.
\end{eqnarray*}
This proves Theorem \ref{thm;finite speed}.

\end{document}